\documentclass{amsart}

\newtheorem{theorem}{\indent Theorem}
\newtheorem{lemma}{\indent Lemma}

\begin{document}

\title[SLLN for $L$-statistics with dependent data]
{On The Strong Law of Large Numbers for $L$-statistics with Dependent Data}

\author{Evgeny Baklanov}

\thanks{The author was supported by
a grant of the President of the Russian Federation for Junior Scientists (Grant MK-2061.2005.1),
the Russian Foundation for Basic Research (Grant 06-01-00738) and INTAS (Grant 03-51-5018).
}

\address{
Department of Mathematics\newline
\phantom{ab} Novosibirsk State University\newline
\phantom{ab} Pirogova st. 2\newline
\phantom{ab} Novosibirsk 630090, Russia}

\email{baklanov@mmf.nsu.ru}

\date{September 27, 2006}

\begin{abstract}
The strong law of large numbers for linear
combinations of functions of order statistics
($L$-statistics) based on weakly dependent random variables is
proven. We also establish the Glivenko--Cantelli theorem for
$\varphi$-mixing sequences of identically distributed random
variables.
\end{abstract}

\maketitle

\section{Introduction}
Let
$X_1$, $X_2$, $\dots$
be a sequence of random variables
with the common distribution function
$F$.
Let us consider the $L$-statistic
$$
L_n=\frac{1}{n}\sum_{i=1}^nc_{ni}h(X_{n:i}),\eqno(1)
$$
where
$X_{n:1}\le\ldots\le X_{n:n}$
are the order statistics based on the sample
$\{X_i, i\le n\}$,
$h$
is a measurable function called a {\it kernel},
$c_{ni}$, $i=1,\dots, n$,
are some constants called {\it weights}.

The aim of this paper is to establish the strong law of large
numbers (SLLN) for $L$-statistics (1) based on sequences of
weakly dependent random variables. The similar problems were
considered in the papers \cite{aaronson} and \cite{gilat}, where the SLLN was proved for
aforementioned $L$-statistics based on stationary ergodic
sequences. For example, in \cite{gilat} the case of linear kernels
($h(x)=x$) and {\it asymptotic regular} weights was considered,
i.~e.
$$
c_{ni}=n\int\limits_{(i-1)/n}^{i/n} J_n(t)\,dt,\eqno(2)
$$
with $J_n$ denoting an integrable function. In addition,
the existence of a function $J$ such that for all
$t\in(0,1)$
$$
\int\limits_0^t J_n(s)\,ds\to\int\limits_0^t J(s)\,ds
$$
was imposed there.
The statistics (1) with linear kernels and {\it regular} weights, i.~e.
$J_n\equiv J$
in (2), were considered in \cite{aaronson}.
In the present paper we relax the regularity assumption on
$c_{ni}$
and, furthermore,  consider the
$L$-statistics (1) based on both stationary ergodic sequences and $\varphi$-mixing
sequences. We also do not impose {\it monotonicity} of the kernel in (1).
Note, that if
$h$
is a monotonic function, then the
$L$-statistic (1) can be represented as a statistic
$$
\frac{1}{n}\sum_{i=1}^nc_{ni}Y_{n:i},
$$
based on a sample
$\{Y_i=h(X_i), i\le n\}$
(see \cite{bakl} for more detail).

As an auxiliary result we obtain the Glivenko--Cantelli theorem
for $\varphi$-mixing sequences.

\section{Notations and Results}
\subsection{Assumptions and notations}
We first introduce our main notations. Let
$F^{-1}(t)=\inf\{x:F(x)\ge t\}$
be the quantile function corresponding to the distribution function
$F$
and let
$U_1$, $U_2$, $\dots$
be a sequence of
uniformly distributed on
$[0,1]$
random variables. Due to the fact
that joint distributions of random vectors
$(X_{n:1},\dots, X_{n:n})$
and
$(F^{-1}(U_{n:1}),\dots, F^{-1}(U_{n:n}))$
coincide, we have that
$$
L_n\stackrel{d}{=}\frac{1}{n}\sum^n_{i=1}c_{ni}H(U_{n:i}),
$$
where
$H(t)=h(F^{-1}(t))$,
and
$\stackrel{d}{=}$
denotes the
equality in distribution. Let us consider a sequence of functions
$c_n(t)=c_{ni}$, $t\in((i-1)/n, i/n]$, $i=1$, \dots, $n$,
$c_n(0)=c_{n1}$.
It is not difficult to see that in this case we have:
$$
L_n=\int\limits_0^1 c_n(t)H(G^{-1}_n(t))\,dt,
$$
where
$G_n^{-1}$
is the quantile function corresponding to the
empirical distribution function
$G_n$
based on the sample
$\{U_i, i\le n\}$.
We also introduce the following notation:
$$
\mu_n=\int\limits_0^1 c_n(t)H(t)\,dt,
$$
$$
C_n(q)=\left\{\begin{array}{ll}n^{-1}\sum\limits_{i=1}^n|c_{ni}|^q\quad&\mbox{if}\quad 1\le q<\infty,\\
\max\limits_{i\le n}|c_{ni}|\quad&\mbox{if}\quad q=\infty.\end{array}\right.
$$
Further we will use the following conditions on the weights
$c_{ni}$
and the function
$H$:\vspace{2mm}

(i) the function
$H$
is continuous on
$[0,1]$
and
$\sup\limits_{n\ge1}C_n(1)<\infty$.

(ii)
${\bf E}|h(X_1)|^p<\infty$
and
$\sup\limits_{n\ge1}C_n(q)<\infty$ ($1\le p<\infty$, $1/p+1/q=1$).
\vspace{2mm}

Assumptions (i) and (ii)
guarantee the existence of
$\mu_n$.
We also note that
$C_n(\infty)=\|c_n\|_{\infty}=\sup\limits_{0\le t\le1}|c_n(t)|$
and
$C_n(q)=\|c_n\|_q^q=\int\limits_0^1 |c_n(t)|^q\,dt$
for
$1\le q<\infty$.

\subsection{SLLN for ergodic and stationary sequences}
Let us formulate our main statement for stationary ergodic sequences.

\begin{theorem}
Let
$\{X_n, n\ge 1\}$
be a strictly stationary and ergodic
sequence and let either {\rm (i)} or {\rm (ii)} hold.
Then, as $n\to\infty$,
$$
L_n-\mu_n\to0\quad\mbox{a. s.}\eqno(3)
$$
\end{theorem}

\textsc{Remark.}
Let us consider the case of regular weights:
$$
c_{ni}=n\int\limits_{(i-1)/n}^{i/n}J(t)\,dt.
$$
Then
$$
L_n=\sum_{i=1}^nH(U_{n:i})\int\limits_{(i-1)/n}^{i/n}J(t)\,dt=
\int\limits_0^1J(t)H(G_n^{-1}(t))\,dt.
$$
Hence, assuming $c_n(t)=J(t)$ in Theorem~1, we have
$$
L_n\to\int\limits_0^1J(t)H(t)\,dt\quad\mbox{a. s.}
$$
Also note that the convergence
$\mu_n\to\mu$, $|\mu|<\infty$,
yields that
$L_n\to\mu$
a.~s. In particular, if
$c_n(t)\to c(t)$
uniformly in
$t\in[0,1]$,
then
$\mu_n\to\int\limits_0^1c(t)H(t)\,dt$.

Without the requirement that the coefficients
$c_{ni}$
are regular
one can easily construct an example when the
assumptions of Theorem~1 are satisfied, but the sequence
$c_{n}(t)$
does not converges in any reasonable sense to a limit function.
Let, for simplicity,
$h(x)=x$
and let
$X_1$ be uniformly distributed on [0, 1].
Set
$c_{ni}=(i-1)\delta_n$
as
$1\le i\le k$
and
$c_{ni}=(2k-i)\delta_n$
as
$k+1\le i\le 2k$, $k=k(n)=[n^{1/2}]$, $\delta_n=n^{-1/2}$.
Thus, the
function $c_{n}(t)$ is defined on the interval $[0, 2k/n]$. On the remaining part of $[0, 1]$ we extend
$c_{n}(t)$ periodically
with period $2k/n$: $c_n(t)=c_n(t-2k/n)$, $2k/n\le t\le 1$
(see also \cite[p.~138]{bakl}).
Note that
$0\le c_{n}(t)\le 1$.
One can show that in
this case
$\mu_n\to 1/4$.
In view of this fact we have that the
assumptions of Theorem~1 are satisfied and, consequently,
$$
L_n\to 1/4\quad\mbox{a. s.}
$$

\subsection{SLLN for $\varphi$-mixing sequences}
We will now formulate our main statement for mixing sequences.
Let us define the mixing coefficients:
$$
\varphi(n)=\sup_{k\ge 1}\sup\{|\mathbf{P}(B|A)-\mathbf{P}(B)|: A\in\mathcal{F}_1^k, B\in\mathcal{F}_{k+n}^\infty, \mathbf{P}(A)>0\},
$$
where
$\mathcal{F}_1^k$
and
$\mathcal{F}_{k+n}^\infty$
denote the
$\sigma$-fields generated by
$\{X_i, 1\le i\le k\}$
and
$\{X_i, i\ge k+n\}$
respectively.
The sequence
$\{X_i, i\ge 1\}$
is called {\it $\varphi$-mixing} (uniform mixing)
if
$\varphi(n)\to0$
as
$n\to\infty$.

\begin{theorem}
Let $\{X_n, n\ge 1\}$ be a $\varphi$-mixing sequence of identically
distributed random variables such that
$$
\sum_{n\ge 1}\varphi^{1/2}(2^n)<\infty,\eqno(4)
$$
and let any of the conditions {\rm (i)} or {\rm (ii)} hold. Then
the statement {\rm (3)} remains true.
\end{theorem}

The proof of Theorem~2 essentially uses the result of the Lemma~1
below. The statement (a) of Lemma~1 is the SLLN for
$\varphi$-mixing sequences. The statement (b) is a
Glivenko--Cantelli-type result for $\varphi$-mixing sequences and
is of independent interest. We note that neither in Theorem~2 nor
in Lemma~1 we do not assume the stationarity of the sequence
$\{X_n\}$.

\begin{lemma}
Let $\{X_n, n\ge 1\}$ be a $\varphi$-mixing sequence of identically
distributed random variables such that the statement {\rm (4)} holds.
Then

(a) for any function
$f$
such that
${\bf E}|f(X_1)|<\infty$,
$$
\frac{1}{n}\sum_{i=1}^nf(X_i)\to\mathbf{E}f(X_1)\quad\mbox{a. s.}\eqno(5)
$$

(b)
$$
\sup_{-\infty<x<\infty}|F_n(x)-F(x)|\to0\quad\mbox{a. s.},\eqno(6)
$$
where
$F_n$ is the empirical distribution function based on the sample
$\{X_i, i\le n\}$.
\end{lemma}

\section{Proofs}
\subsection{Proof of Theorem 1}
%Now let us proceed to the proof of Theorem~1.

\begin{lemma}
Let the function
$H$
be continuous on
$[0,1]$.
Then
$$
\sup_{0\le t\le1}|H(G_n^{-1}(t))-H(t)|\to0\quad\mbox{a. s.}\eqno(7)
$$
\end{lemma}

\textsc{Proof} of Lemma~2.
Using the equality
$$
\sup_{0\le t\le1}|G^{-1}_n(t)-t|=\sup_{0\le t\le1}|G_n(t)-t|
$$
(see, for example, \cite[p.~95]{sh-well})
and the Glivenko--Cantelli theorem for stationary ergodic sequences, we get
$$
\sup_{0\le t\le1}|G^{-1}_n(t)-t|\to0\quad\mbox{a. s.},
$$
i.~e.
$G^{-1}_n(t)\to t$
a.~s. uniformly in
$t\in[0,1]$
as
$n\to\infty$.
Since the function
$H$
is uniformly continuous on the compact
$[0,1]$,
it follows that
$H(G_n^{-1}(t))\to H(t)$
a. s. uniformly in
$t\in[0,1]$.
This concludes the proof.

Let the condition (i) hold. Now, by Lemma~2,
$$
|L_n-\mu_n|\le\int\limits_0^1|c_n(t)||H(G_n^{-1}(t))-H(t)|\,dt
$$
$$
\le
C_n(1)\sup_{0\le t\le 1}|H(G_n^{-1}(t))-H(t)|\to0\quad\mbox{a. s.}
$$
Consequently, the proof of Theorem~1 for the first case is complete.

\begin{lemma}
Let
${\bf E}|h(X_1)|^p<\infty$.
Then
$$
\int\limits_0^1|H(G_n^{-1}(t))-H(t)|^p\,dt\to0\quad\mbox{a. s.}\eqno(8)
$$
\end{lemma}

\textsc{Proof} of Lemma~3.
First note that the set of all continuous on the interval
$[0,1]$
functions is everywhere dense in
$L_p[0,1]$, $1\le p<\infty$.
Therefore, for any
$\varepsilon>0$
and any function
$f\in L_p[0,1]$
there exists a continuous on
$[0,1]$
function
$f_{\varepsilon}$
such that
$\int\limits_0^1|f(t)-f_{\varepsilon}(t)|^p\,dt<\varepsilon$.
Since
$\mathbf{E}|h(X_1)|^p=\int\limits_0^1|H(t)|^p\,dt<\infty$,
this implies that there exists a continuous on
$[0,1]$
function
$H_{\varepsilon}$
such that
$$
\int\limits_0^1|H(t)-H_{\varepsilon}(t)|^pdt<\varepsilon/2.
$$
Further,
$$
\int\limits_0^1|H(G_n^{-1}(t))-H(t)|^p\,dt\le3^{p-1}
\int\limits_0^1|H(t)-H_{\varepsilon}(t)|^p\,dt
$$
$$
+3^{p-1}\int\limits_0^1|H(G_n^{-1}(t))-H_{\varepsilon}(G_n^{-1}(t))|^p\,dt+
3^{p-1}\int\limits_0^1|H_{\varepsilon}(G_n^{-1}(t))-H_{\varepsilon}(t)|^p\,dt.\eqno(9)
$$
From Lemma~2 it follows that
$H_{\varepsilon}(G_n^{-1}(t))\to H_{\varepsilon}(t)$
a.~s. uniformly in
$t$
as
$n\to\infty$.
Hence, the last integral on the right hand side of (9) converges to zero a. s. as
$n\to\infty$.
Now let us consider the second integral. By ergodic theorem for stationary sequences,
\begin{eqnarray*}
&&
\int\limits_0^1|H(G_n^{-1}(t))-H_{\varepsilon}(G_n^{-1}(t))|^p\,dt\\
&=&\frac{1}{n}\sum_{i=1}^{n}|H(U_i)-H_{\varepsilon}(U_i)|^p
\to_{\mbox{a. s.}}\mathbf{E}|H(U_1)-H_{\varepsilon}(U_1)|^p\\
&=&\int\limits_0^1|H(t)-H_{\varepsilon}(t)|^p\,dt<\varepsilon/2.
\end{eqnarray*}
Consequently,
$$
\limsup_{n\to\infty}\int\limits_0^1|H(G_n^{-1}(t))-H(t)|\,dt<3^{p-1}\varepsilon\quad\mbox{a. s.}
$$
Since
$\varepsilon$
is arbitrary, we obtain (8).

Now let the assumption (ii) hold. Using H\"{o}lder's inequality, we get
$$
|L_n-\mu_n|\le C_n^{1/q}(q)\left(\int\limits_0^1|H(G_n^{-1}(t)-H(t))|^p\,dt\right)^{1/p}\quad\mbox{for }p>1,
$$
and
$$
|L_n-\mu_n|\le C_n(\infty)\int\limits_0^1|H(G_n^{-1}(t)-H(t))|\,dt\quad\mbox{for }p=1.
$$
The statement (3) follows from Lemma~3.
This completes the proof of Theorem~1.

\subsection{Proof of Theorem 2}
We now prove Lemma~1.
Note that for any measurable function $f$
the sequence
$\{f(X_n), n\ge 1\}$
has its $\varphi$-mixing coefficient bounded by the corresponding
coefficient of the initial sequence, since
for any measurable $f$
the $\sigma$-field generated by $\{f(X_n), n\ge 1\}$ is contained in the
$\sigma$-field generated by $\{X_n, n\ge 1\}$.
Therefore, if the sequence $\{X_n, n\ge 1\}$ is $\varphi$-mixing,
then so is the sequence $\{f(X_n), n\ge 1\}$.
Hence, the condition (4) holds
for mixing coefficients of the sequence
$\{f(X_n), n\ge 1\}$.
The statement (5) follows from the SLLN for
$\varphi$-mixing sequences
(see \cite[p.~200]{linlu}).

The statement (6) is an immediate corollary of (5) and classical Glivenko--Cantelli theorem.

The proof of Theorem~2 is similar to the proof of Theorem~1. Indeed, the statement (7)
follows from the Glivenko--Cantelli theorem (6); using the SLLN (5), we get the statement (8).
Thus the proof of Theorem~2 is complete.

\end{document}